\documentclass{article}

\usepackage{arxiv}

\usepackage[utf8]{inputenc}
\usepackage[T1]{fontenc}
\usepackage{amsmath}
\usepackage{amssymb}
\usepackage{booktabs}
\usepackage{graphicx}
\usepackage{xcolor}
\usepackage{microtype}
\usepackage[numbers]{natbib}
\usepackage{doi}
\usepackage{hyperref}
\usepackage{url}

\title{Fast Karhunen--Loève Expansions via \\
  FFT-Accelerated Toeplitz Operators}

\date{\today}

\usepackage{authblk}

\author[1]{Nils Wildt}
\author[1]{Wolfgang Nowak}
\affil[1]{Department of Stochastic Simulation and Safety Research for Hydrosystems, University of Stuttgart}

\renewcommand{\shorttitle}{Fast Karhunen--Loève Expansions via FFT-Accelerated Toeplitz Operators}

\definecolor{linkblue}{HTML}{1C3F6E}
\hypersetup{
  colorlinks=true, linkcolor=linkblue, citecolor=linkblue, urlcolor=linkblue,
  pdftitle={Fast Karhunen-Loève Expansions via FFT-Accelerated Toeplitz Operators},
  pdfauthor={Nils Wildt, Wolfgang Nowak},
}

\begin{document}
\maketitle

\begin{abstract}
Gaussian random fields are a versatile tool used in
the fields of stochastic PDEs, uncertainty quantification, and
geostatistical simulation. One way to obtain them is to use a truncated
Karhunen--Loève expansion (KLE). Computing the expansion requires the
leading eigenpairs of an \(N \times N\) covariance matrix, where \(N\)
is the total number of grid cells. These are usually computed with a
Krylov eigensolver, which relies on the covariance operator only within
matrix--vector products. Stored densely, the matrix takes
\(\mathcal{O}(N^{2})\) memory and each product \(\mathcal{O}(N^{2})\)
time. For a stationary kernel on an equispaced grid, the covariance
matrix becomes (block-) Toeplitz and the product evaluates in
\(\mathcal{O}(N\log N)\) time using FFT-based circulant embedding,
without the need to assemble the dense matrix. In a matched
single-threaded comparison, the median speedup of the eigensolve grows
from \(18 \times\) at \(N = 4096\) to \(183 \times\) at \(N = 2^{15}\).
This makes it possible to compute discretized fields that would
otherwise be infeasible to compute in the standard formulation. We show
that the same construction carries over to non-separable kernels as well
as \(d\) dimensions, using block-Toeplitz matrices. We extend it to
piecewise-constant fields on arbitrary domains given as subsets of a
tensor grid. Computational savings grow with problem size, and storage
drops from \(\mathcal{O}(N^{2})\) to \(\mathcal{O}(2^{d}N)\).
\end{abstract}

\section{Motivation}

A truncated KLE expands a Gaussian random field into eigenfunctions
weighted by independent standard-normal coefficients and underlies
random-field generation for stochastic PDEs, polynomial-chaos
uncertainty quantification, and geostatistical simulation
\citep{lord2014, ghanem1991}. Direct circulant-embedding samplers
generate individual realizations without solving an eigenproblem
\citep{dietrich1997, graham2018}. The KLE is used when the eigenpairs
themselves are needed, e.g.~for polynomial-chaos expansions, and
whenever one truncated basis is reused across many realizations or
parameter values. We consider fields that are piecewise constant on the
\(N\) equal-volume cells of an equispaced tensor grid, so a field is a
random vector of \(N\) values, one per cell, sampled at the cell's grid
point. Its covariance is the \(N \times N\) matrix \(K\) of kernel
values between those points. For a finite random vector the KLE is
simply the eigendecomposition of its covariance matrix. All cells have
the same volume, so no quadrature weights or mass matrix enter. The
standard eigenproblem \(Ku = \lambda u\) is the discrete KLE itself.
Computing the truncated expansion means finding the leading \(m\)
eigenpairs of \(K\), where \(m\) is the truncation rank.

Once \(N\) exceeds a few thousand, Krylov methods (Lanczos for a
symmetric operator, implicitly restarted Arnoldi in general) are
commonly used rather than a full eigendecomposition. Such solvers touch
\(K\) only through the products \(v \mapsto Kv\). Most standard
implementation work this way e.g.~ARPACK \citep{arpack}. The cost of the
solve is the cost of that one operation. With \(K\) stored densely, each
product costs \(\mathcal{O}(N^{2})\) time and \(K\) itself
\(\mathcal{O}(N^{2})\) memory, however efficient the eigensolver.

For a stationary kernel on an equispaced grid,
\(K_{ij} = c\left( x_{i} - x_{j} \right)\) depends only on the
displacement between cells, so \(K\) becomes Toeplitz. Its
matrix--vector product reduces to \(\mathcal{O}(N\log N)\) by embedding
\(K\) in a circulant matrix and diagonalizing the embedding with the FFT
\citep{dietrich1997}. A Krylov solver only uses \(K\) in this
multiplication, so the FFT replacement for the dense matrix--vector
multiplication in any implementation leaves its convergence behavior and
guarantees untouched (see
Section~\ref{sec-discussion}).

We believe that the combination of the FFT/Toeplitz operator inside a
Krylov eigensolver, for truncated KLE of stationary random fields, might
well have been used by practitioners, as the idea is straightforward. We
are, however, not aware of published work on this specific combination.
Khoromskij, Litvinenko \& Matthies \citep{khoromskij2009} noted that FFT
techniques apply to the resulting block-Toeplitz covariance matrices on
uniform rectangular grids. They pursued hierarchical-matrix products for
more general settings instead. Here we work out the FFT-based route in
full, for arbitrary fixed spatial dimension. The closest precedent for
the FFT/Toeplitz idea on the same class of matrices is due to Fritz,
Neuweiler \& Nowak \citep{fritz2009}. They built an FFT-based Toeplitz
solver for universal kriging systems, a linear solve rather than an
eigenproblem. Beyond hierarchical matrices, \citet{schwab2006}
accelerate the KLE eigensolve itself with a generalized fast multipole
method, reaching log-linear cost for analytic covariances on general
polyhedral domains. Our operator targets the narrower case of a
stationary kernel on an equispaced tensor grid, where the FFT gives an
exact matrix--vector product with no multipole machinery. Randomized
matrix-free eigensolvers consume the operator only within matrix--vector
products and thus inherit the same \(\mathcal{O}(N\log N)\) cost. This
holds for subspace iteration built for this very class of problem
\citep{saibaba2016} as well as for randomized low-rank algorithms in
general \citep{halko2011}.

\section{Method}

\label{sec-method}

A Krylov eigensolver never needs the full and explicit covariance
matrix. For a stationary kernel on an equispaced grid it relies only on
matrix--vector products, and each of them can be evaluated in
\(\mathcal{O}(N\log N)\) time.

\textbf{The substitution:} Starting from \(v_{0}\), such a solver builds
everything it needs from the sequence
\(v_{0},Kv_{0},K^{2}v_{0},\ldots\). It only ever applies \(K\) to
vectors and never inspects the matrix entries. In exact arithmetic, two
representations computing the same products, \(A_{1}v = A_{2}v = Kv\)
for every \(v\), produce identical iterates, Ritz values, and eigenpairs
from the same \(v_{0}\). Only the cost of the product differs. The dense
matrix can therefore be replaced by any representation of the same
operator without changing the solver. In floating point arithmetic, the
two evaluations round differently and have different iteration
histories, but the eigenproblem and its convergence theory are unchanged
(Section~\ref{sec-discussion}).

\textbf{The Toeplitz representation:} Embed the length-\(N\) Toeplitz
generator \(c\) in a circulant vector \(c_{\mathrm{circ}}\) of length
\(M \geq 2N - 1\) (padded to a fast FFT length), and zero-pad \(v\) to
the same length. Then we have
\[Kv = \mathrm{IFFT}\left( \,\mathrm{FFT}\left( c_{\mathrm{circ}} \right) \cdot \mathrm{FFT}\left( v_{\mathrm{padded}} \right)\, \right)_{1:N},\]
where the product is evaluated at
\(\mathcal{O}(M\log M) = \mathcal{O}(N\log N)\) cost and
\(\mathcal{O}(M)\) working storage, against \(\mathcal{O}(N^{2})\) for a
dense \(K\). Circulant-embedding samplers need a positive semi-definite
embedding and may have to enlarge the padding until it is
\citep{dietrich1997, graham2018}. Here only the product is required. The
identity holds for every \(M \geq 2N - 1\), whatever the embedding's
spectrum. The saving has a simple source. \(K\) has only \(2N - 1\)
distinct entries, one per displacement. The dense matrix stores each of
them about \(N\) times and the product reads all \(N^{2}\) copies. The
structured product stores them once, as
\(\text{FFT}\left( c_{\mathrm{circ}} \right)\), and touches
\(\mathcal{O}(M)\) numbers per call. The arithmetic falls as well. The
product is a convolution, and in the Fourier basis a convolution is
\(M\) pointwise multiplications instead of \(N^{2}\) cross terms. The
FFT provides that change of basis in \(\mathcal{O}(M\log M)\)
\citep{golub2013}.

\textbf{The \(d\)-dimensional generalization:} The same construction
applies one axis at a time. Let the tensor grid have \(n_{k}\) cells
along axis \(k\) with spacing \(h_{k}\), so that \(N = \prod_{k}n_{k}\)
cells sit at the positions
\(x_{i} = \left( i_{1}h_{1},\ldots,i_{d}h_{d} \right)\). Stationarity
makes the covariance depend only on the displacement,
\(K_{ij} = c\left( x_{i} - x_{j} \right)\). With the cells ordered axis
by axis, \(K\) is therefore Toeplitz in the index of each axis. It is an
\(n_{d} \times n_{d}\) Toeplitz arrangement of blocks, each block again
Toeplitz one axis down, and so on through all \(d\) axes. This is a
\(d\)-level Toeplitz matrix \citep{golub2013}. For \(d = 2\) this is the
familiar block-Toeplitz with Toeplitz blocks (BTTB), and we use that
name for every \(d\).

Applying \(K\) needs a circulant array whose leading
\(n_{1} \times \ldots \times n_{d}\) block is \(K\). We pad axis \(k\)
to a length \(M_{k} \geq 2n_{k} - 1\) (in practice the next 2-3-5-smooth
length, so that the transforms stay fast) and read each padded index
\(j_{k} \in \left\{ 0,\ldots,M_{k} - 1 \right\}\) as a signed
displacement modulo \(M_{k}\):
\[\delta_{k}\left( j_{k} \right) = \begin{cases}
j_{k} & \text{if }0 \leq j_{k} \leq n_{k} - 1, \\
j_{k} - M_{k} & \text{if }M_{k} - n_{k} + 1 \leq j_{k} \leq M_{k} - 1, \\
\text{undefined} & \text{otherwise.}
\end{cases}\] The two defined ranges together give all \(2n_{k} - 1\)
displacements \(- \left( n_{k} - 1 \right),\ldots,n_{k} - 1\) between
two cells on axis \(k\). The remaining \(M_{k} - 2n_{k} + 1\) indices
correspond to no cell pair. There are none when \(M_{k} = 2n_{k} - 1\).
The padded kernel array
\(\Lambda_{\mathrm{pad}} \in {\mathbb{R}}^{M_{1} \times \ldots \times M_{d}}\)
is
\[\Lambda_{\mathrm{pad}}\left\lbrack j_{1},\ldots,j_{d} \right\rbrack = \begin{cases}
c\left( \delta_{1}\left( j_{1} \right)h_{1},\ldots,\delta_{d}\left( j_{d} \right)h_{d} \right) & \text{if every }\delta_{k}\left( j_{k} \right)\text{ is defined,} \\
0 & \text{otherwise.}
\end{cases}\] Because \(c\) is even, the entries at negative
displacements are those at positive displacements reflected, so a single
rule serves every axis and every dimension. With \(M = \prod_{k}M_{k}\),
we can precompute its transform
\(\Lambda = \mathrm{FFT}\left( \Lambda_{\mathrm{pad}} \right)\).

The \(d\)-dimensional DFT diagonalizes the circulant generated by
\(\Lambda_{\mathrm{pad}}\): zero-padding \(v \in {\mathbb{R}}^{N}\) to
\(v_{\mathrm{pad}} \in {\mathbb{R}}^{M_{1} \times \ldots \times M_{d}}\),
with the entries of \(v\) at the indices \(0 \leq j_{k} \leq n_{k} - 1\)
and zeros elsewhere,
\[Kv = \,\left( \mathrm{IFFT}\left( \Lambda \cdot \mathrm{FFT}\left( v_{\mathrm{pad}} \right) \right)\, \right)_{\left\lbrack 0:n_{1} - 1,\ldots,0:n_{d} - 1 \right\rbrack},\]
with \(\Lambda\) taken unnormalized and the inverse transform carrying
the \(1/M\), as in the one-dimensional case. The identity is exact
because the leading \(n_{1} \times \ldots \times n_{d}\) block of that
circulant is \(K\). By construction, \(\Lambda_{\mathrm{pad}}\) at the
displacement between two cells is their kernel value, so restricting the
padded product to those indices returns \(Kv\) unchanged. Since \(c\) is
even, \(\Lambda_{\mathrm{pad}}\) is symmetric about the origin and
\(\Lambda\) is real, but it may be negative. The cost is two
\(d\)-dimensional transforms of size \(M\) plus \(M\) pointwise
multiplications, \(\mathcal{O}(M\log M)\) against \(\mathcal{O}(N^{2})\)
for a dense product, with \(\mathcal{O}(M)\) memory required. Since
padding roughly doubles each axis, \(M_{k} \approx 2n_{k}\), the padded
size is \(M = \prod_{k}M_{k} \approx 2^{d}\prod_{k}n_{k} = 2^{d}N\). For
fixed \(d\), \(M\) is a constant multiple of \(N\), and \(2^{d}\) is the
price of the embedding.

\textbf{Masked domains:} The same operator covers domains that are not
boxes. Let
\(\Omega \subseteq \left\{ 1,\ldots,N_{\mathrm{full}} \right\}\) index the
active cells of a bounding-box tensor grid with \(N_{\mathrm{full}}\)
cells, and let \(v\) be defined on \(\Omega\). Write \(E\) for
zero-extension from \(\Omega\) to the box and \(R = E^{\top}\) for
restriction back to \(\Omega\). Then at every active cell \(i\)
\[(RKEv)_{i} = \sum_{j \in \Omega}K_{ij}v_{j} = \left( K_{\Omega\Omega}v \right)_{i},\]
since \(Ev\) vanishes off \(\Omega\). Hence \(RKE = K_{\Omega\Omega}\)
for any matrix \(K\). Stationarity is needed only for the cost. \(K\) on
the bounding box is BTTB, so \(RKEv\) is the full-box product of the
previous paragraph, applied to the zero-extended vector and read off on
\(\Omega\). The FFTs of size \(M\) and the working storage are those of
the bounding box and independent of \(\Omega\). A mask covering any
fraction of its box costs the same as the full box. Transposing \(RKE\)
swaps the restriction and the extension, and these are transposes of
each other. Therfore \(RKE\) is symmetric whenever \(K\) is and the
solver still solves on a symmetric operator on the masked domain.

\section{Results}

\label{sec-results}

We compare three ways of computing the leading \(m\) eigenpairs of
\(K\). \textbf{Full dense} assembles \(K\) and calls LAPACK's
\texttt{eigen} for all \(N\) eigenpairs. It serves as the reference and,
at \(\mathcal{O}(N^{3})\), is not competitive. \textbf{Dense Krylov}
assembles \(K\) and hands it to ARPACK's \texttt{eigs}. This is what
standard KLE implementations do (e.g.~in the widely used
GaussianRandomFields.jl Julia package \citep{grfjl}) and is our
baseline. \textbf{BTTB Krylov} hands the matrix-free operator of
Section~\ref{sec-method} to the same \texttt{eigs} call
with the same settings, so the two Krylov variants differ only in how
\(Kv\) is evaluated.

\begin{figure}[t]
\centering
\includegraphics[width=\linewidth]{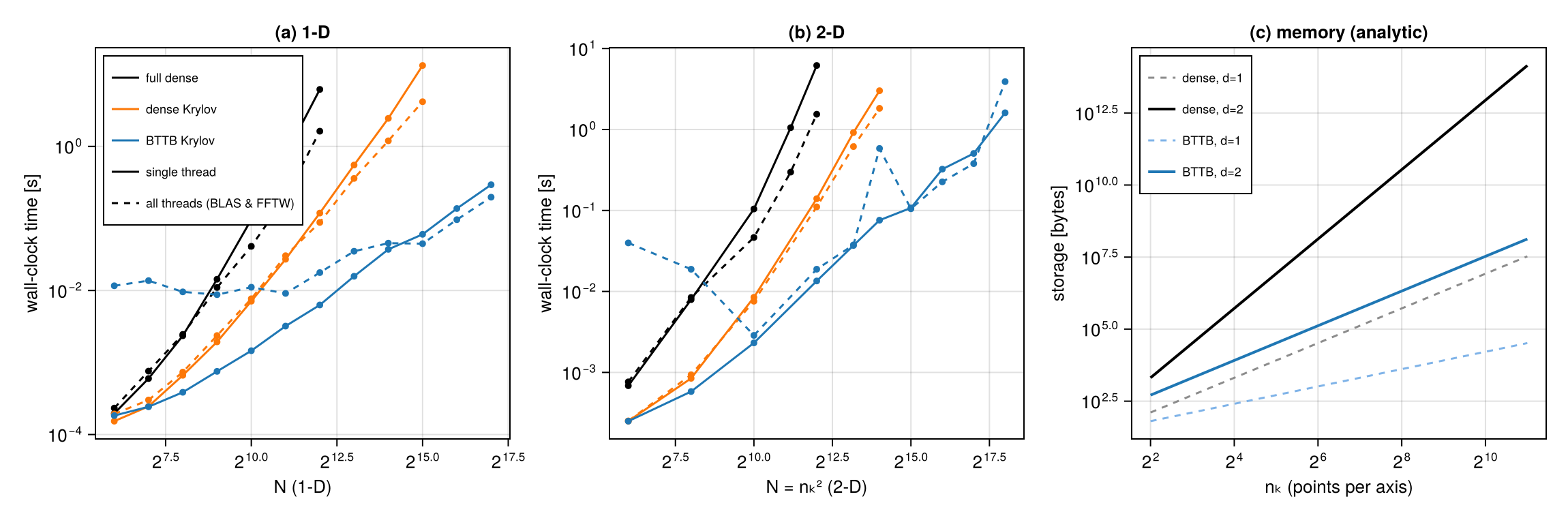}
\caption{(a) Wall-clock time of full dense, dense Krylov, and BTTB
Krylov to solve for \(m = 20\) eigenpairs vs. \(N\), 1-D. (b) Same, 2-D,
\(N = n_{k}^{2}\). In (a) and (b), solid lines are single-threaded,
dashed lines use all ten threads for BLAS and FFTW alike. (c) Analytic
storage vs. per-axis resolution \(n_{k}\). Dense storage grows as
\(n_{k}^{2d}\), the BTTB generator as \(n_{k}^{d}\), so the storage gap
widens with \(d\) at fixed per-axis resolution.}
\label{fig-benchmark}
\end{figure}

We validated the implementation against full dense on 1-D, 2-D, and 3-D
grids, on a non-separable kernel \(C(h) = \exp( - h^{\top}Ah)\) with
\(A = \begin{pmatrix}
18 & 10 \\
10 & 18
\end{pmatrix}\), whose precision matrix couples the axes, and on
L-shaped and circular-cutout masked domains
(Table~\ref{tbl-validation}). In every case BTTB
Krylov reproduces the full dense eigenvalues to a relative error of
approximately machine epsilon and the leading subspace to a principal
angle below \(6 \times 10^{- 8}\) rad, approximately the accuracy the
eigensolver itself can obtain. The rest of this section reports timings
(Figure~\ref{fig-benchmark}) and the masked-domain
experiments.

\begin{table}[t]
\centering
\small
\begin{tabular}{lrrrrr}
\toprule
\textbf{Configuration} & \(N\) & \(m\) & \textbf{eigenvalue} &
\textbf{principal angle} & \textbf{residual} \\
\midrule
1-D grid & \(256\) & \(15\) & \(1.5 \cdot 10^{- 15}\) &
\(2.6 \cdot 10^{- 8}\) & \(7.6 \cdot 10^{- 16}\) \\
1-D grid & \(4096\) & \(15\) & \(1.3 \cdot 10^{- 15}\) &
\(3.9 \cdot 10^{- 8}\) & \(1.1 \cdot 10^{- 15}\) \\
2-D grid & \(1024\) & \(15\) & \(1.4 \cdot 10^{- 15}\) &
\(3.0 \cdot 10^{- 8}\) & \(1.2 \cdot 10^{- 15}\) \\
3-D grid & \(512\) & \(10\) & \(3.1 \cdot 10^{- 15}\) &
\(5.6 \cdot 10^{- 8}\) & \(1.8 \cdot 10^{- 15}\) \\
2-D non-separable & \(400\) & \(12\) & \(1.6 \cdot 10^{- 15}\) & --- &
--- \\
2-D L-shape mask & \(1728\) & \(6\) & \(9.5 \cdot 10^{- 16}\) & --- &
\(1.6 \cdot 10^{- 15}\) \\
2-D circular-cutout mask & \(1856\) & \(6\) & \(1.5 \cdot 10^{- 15}\) &
--- & \(1.8 \cdot 10^{- 15}\) \\ \\
\bottomrule
\end{tabular}
\caption{BTTB Krylov against full dense. Eigenvalue:
\(\max\limits_{i}|\lambda_{i}^{\mathrm{BTTB}} - \lambda_{i}^{\mathrm{full}}|/\lambda_{i}^{\mathrm{full}}\)
for \(i = 1,\ldots,m\). Principal angle: largest \(\arccos\) of the
singular values of \(Q_{\mathrm{full}}^{\top}Q_{\mathrm{BTTB}}\). Residual:
\(\max\limits_{k}\left\| {Kv_{k} - \lambda_{k}v_{k}} \right\|/\lambda_{1}\).
Grids use the squared-exponential kernel (\(\ell = 0.05\) in 1-D,
\(0.2\) otherwise), masked domains \(\ell = 0.15\).}
\label{tbl-validation}
\end{table}

\textbf{Performance:} In
Figure~\ref{fig-benchmark}(a) we sweep \(N\) in 1-D
for a squared-exponential kernel (\(\ell = 0.05\)) and \(m = 20\)
eigenpairs.\footnote{Apple M4 Pro, macOS 26.6.2 (arm64); Julia 1.13.0,
  \texttt{Float64} precision. Every point was measured with one thread
  and all ten threads for FFTW as well as BLAS. Timings are single timed
  calls after an untimed warm-up that excludes FFTW planning and JIT
  compilation. We report the medians of five paired repetitions.} The
covariance construction sits outside the timed region for both methods.
This favors dense Krylov, as assembling \(K\) costs
\(\mathcal{O}(N^{2})\) kernel evaluations against \(\mathcal{O}(N)\) for
the Toeplitz generator. From the paired medians, the single-threaded
speedup is \(18 \times\) at \(N = 4096\), \(36 \times\) at \(8192\),
\(66 \times\) at \(16384\), and \(183 \times\) at \(N = 2^{15}\).
Threading reduces the difference only by constant factors. At
\(N = 2^{15}\), ten BLAS threads give dense Krylov a \(1.8 \times\)
speedup. BTTB Krylov gains almost nothing, and threaded FFTW can lose on
small transforms. The all-threads pairing yields \(4 \times\),
\(12 \times\), \(30 \times\), and \(107 \times\) speedup, but BTTB
Krylov is fastest when run single-threaded. The single-thread pairing is
fair. Below \(N \approx 2^{8}\) the two methods are comparable, and
above it the difference increases.

All speedups compare dense Krylov with BTTB Krylov. Full dense appears
in panel (a) only to show that it is slower, as expected. Dense Krylov
stops at \(N = 2^{15}\), because at \(N = 2^{17}\) the covariance matrix
would need 128 GiB of memory versus 1 MiB for its Toeplitz generator.
The \(N = 2^{16},2^{17}\) BTTB Krylov points demonstrate feasibility.
The 2-D sweep in Figure~\ref{fig-benchmark}(b)
(\(\ell = 0.2\)) repeats the pattern up to \(N = 2^{14} = 16384\).
There, the single-threaded speedup is about \(40 \times\). Paired sweeps
beyond 2-D give the same picture, with \(12 \times\) at \(N = 8000\) in
3-D and \(24 \times\) at \(N = 9216\) for the non-separable kernel,
single-threaded (\(6 \times\) and \(12 \times\) with all threads). A
comparably wide 3-D dense sweep would be prohibitive. The comparison
continues with BTTB Krylov alone to \(N \approx 2.6 \times 10^{5}\),
where a dense \(K\) would need 512 GiB. At fixed per-axis resolution
\(n_{k}\), dense storage scales as \(n_{k}^{2d}\) versus \(n_{k}^{d}\)
for the BTTB generator (about \(2^{d}n_{k}^{d}\) with embedding, see
Figure~\ref{fig-benchmark}(c)).

\begin{table}[t]
\centering
\begin{tabular}{lccc}
\toprule
\textbf{Method} & \textbf{Covariance storage} & \textbf{Cost per
matrix--vector product} & \textbf{Operator work in solve} \\
\midrule
Full dense & \(\mathcal{O}(N^{2})\) & --- & \(\mathcal{O}(N^{3})\)
overall \\
Dense Krylov & \(\mathcal{O}(N^{2})\) & \(\mathcal{O}(N^{2})\) &
\(\mathcal{O}(kN^{2})\) \\
BTTB Krylov & \(\mathcal{O}(M)\) & \(\mathcal{O}(M\log M)\) &
\(\mathcal{O}(kM\log M)\) \\ \\
\bottomrule
\end{tabular}
\caption{Asymptotic covariance-operator costs. \(k\) is the number of
matrix--vector products performed by the Krylov eigensolver. \(M\) is
the circulant-embedding size, \(M \approx 2^{d}N\) for fixed \(d\) when
each axis is approximately doubled. Solver-internal orthogonalization
and projected eigenproblems are omitted because they occur identically
for both Krylov methods.}
\label{tbl-cost}
\end{table}

Shorter correlation lengths (rougher fields) need a larger \(m\) for a
given truncation error, hence more restart cycles and a larger count
\(k\) of matrix--vector products. But \(k\) multiplies the per-product
cost of both operators identically
(Table~\ref{tbl-cost}), so the speedup ratio is, to
leading order, independent of \(k\) and of correlation length. Therefore
we expect a similar order of speedup for other field configurations.

\textbf{Irregular domains:} By Section~\ref{sec-method},
restricting to a masked domain changes nothing about how the operator
works. The FFTs run on the bounding box, and the mask only restricts the
input into the box and extends the output from it. So there is no new
speed or memory change to be expected. We show that the restriction and
extension are implemented correctly and that the cost is indeed that of
the bounding box. First, we compared BTTB Krylov against full dense on
\(K_{\Omega\Omega}\), assembled by direct pairwise evaluation, on an
L-shaped domain (\(48 \times 48\) grid, one quadrant removed,
\(N_{\Omega} = 1728\), \(75\%\) occupancy) and on a square with a
circular cutout (\(N_{\Omega} = 1856\), \(81\%\) occupancy);
Table~\ref{tbl-validation} lists the eigenvalue and
residual agreement, and the same checks pass on the unit square. On
masked domains we report residuals rather than principal angles, because
the principal angle presumes that the top-\(m\) subspace is unique. A
symmetric shape has exactly repeated eigenvalues, and such a pair can
straddle the truncation. In the validation runs of
Table~\ref{tbl-validation}, \(m = 6\) and the
circular cutout has \(\lambda_{6} = \lambda_{7}\). Any vector in that
two-dimensional eigenspace is then a legitimate sixth eigenvector. Full
dense and BTTB Krylov may return different ones, and a subspace
comparison would falsely appear as an error. The residual
\(\left\| {Kv_{k} - \lambda_{k}v_{k}} \right\|/\lambda_{1}\) compares
each returned pair against the operator directly and is therefore
unaffected.

\begin{figure}[t]
\centering
\includegraphics[width=\linewidth]{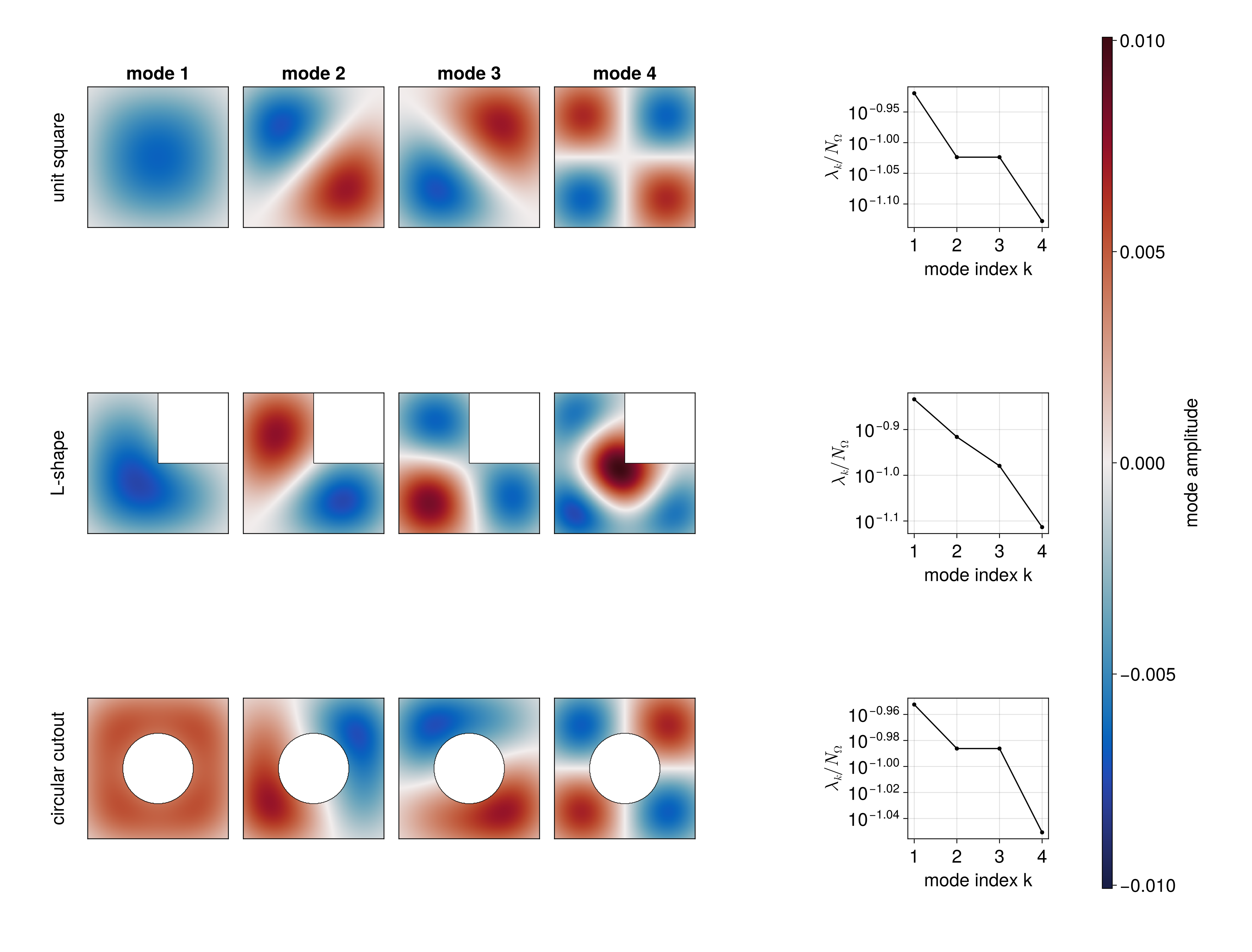}
\caption{Leading four eigenmodes of the squared-exponential KLE
(\(\ell = 0.15\)) on various shapes. Each row ends with the
\(N_{\Omega}\)-normalized eigenvalue decay. All mode panels share one
symmetric color scale.}
\label{fig-modes-standard}
\end{figure}

\begin{figure}[t]
\centering
\includegraphics[width=\linewidth]{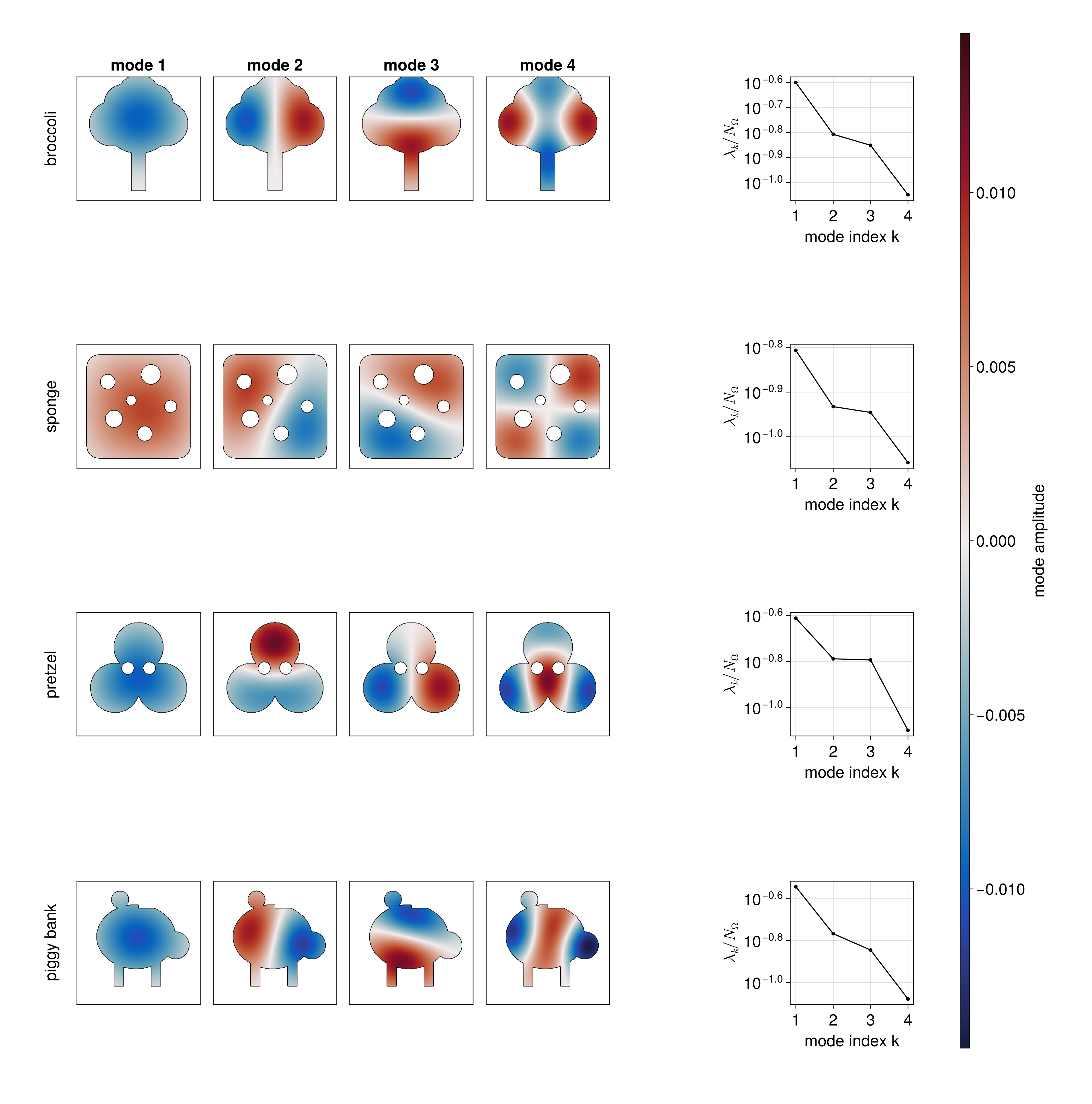}
\caption{The same on broccoli, sponge, pretzel, and piggy-bank shapes,
in the layout of
Figure~\ref{fig-modes-standard}.}
\label{fig-modes-gallery}
\end{figure}

Second, we timed the L-shape against its bounding box and against dense
Krylov on the assembled \(K_{\Omega\Omega}\). The masked matrix--vector
product costs as much as the the full-box costs, and the masked solve
also takes the same number of matrix--vector products as dense Krylov at
every size agreeing with it to machine precision. With the operator
validated and its cost confirmed, the eigenmodes on arbitrary shapes
come for free. Just to please the eye,
Figure~\ref{fig-modes-standard} shows the
leading four on the unit square, the L-shape, and the circular cutout.
Figure~\ref{fig-modes-gallery} adds some
broccoli, a sponge, a pretzel, and a piggy bank. Each has its own mask
but is computed by the same construction.

\section{Discussion and Conclusion}

\label{sec-discussion}

The contribution is deliberately narrow. We expose the Toeplitz/BTTB
matrix multiplication as a matrix-free operator to a standard Krylov
eigensolver, wherever the covariance has that structure. We implement it
generically in the spatial dimension \(d\), and extend it to
piecewise-constant fields on masked domains. The numerical study
confirms agreement to machine precision and storage reduced from
\(\mathcal{O}(N^{2})\) to \(\mathcal{O}(2^{d}N)\). The gain comes from
the redundancy of \(K\) exploited by evaluating the product as a
convolution in the Fourier basis
(Section~\ref{sec-method}).

In the 1-D squared-exponential experiment the median single-threaded
speedup rises from \(18 \times\) at \(N = 4096\) to \(183 \times\) at
\(N = 2^{15}\). The same operator gives \(24 \times\) at \(N = 9216\)
for a non-separable stationary kernel. On masked domains the product
costs the same as the full-box product, and the solve comes in at about
three quarters of the full-box solve, since it converges in fewer
products.

We tested the exact-arithmetic prediction of
Section~\ref{sec-method} directly, with ten paired solves
per size in each dimension. Both operators started from the same seeded
\(v_{0}\) with identical solver settings, so only the floating-point
evaluation of \(Kv\) differs between them. In 1-D and 2-D they needed
the same number of matrix--vector products at every size and in every
run, so the solve inherits the full speedup of the product. In 3-D the
counts differ, with BTTB about 9\% higher on average. That trims a
comparable fraction from the solve speedup but leaves the \(12 \times\)
of Section~\ref{sec-results} essentially intact. The
usual break-even at small \(N\) remains.

The substitution is available exactly when the covariance is
(block-)Toeplitz, that is, for a stationary kernel on an equispaced
tensor grid. Non-stationary kernels, meshes that are unstructured or
locally refined, and grids equispaced only along some axes all fall
outside the scope of this implementation, because \(K\) would no longer
be a BTTB matrix.

Beyond this regime, hierarchical-matrix \citep{khoromskij2009} and
fast-multipole \citep{schwab2006} operators remain the tools of choice.
Within it, randomized matrix-free eigensolvers
\citep{saibaba2016, halko2011} are complementary rather than competing,
since they too touch the covariance only through matrix--vector
products. Any solver with that property inherits the same
\(\mathcal{O}(N\log N)\) cost.

\subsection{Data Availability Statement}

Code, raw benchmark data, results, and the figure-generation scripts
accompany this letter as a self-contained bundle\footnote{\url{https://doi.org/10.18419/DARUS-6487}}.

\section*{AI Disclosure Statement}

The research questions, underlying concept and solutions are the
authors' own, developed through discussions and research, and older than
modern LLMs. However, generative AI tools (Anthropic's Opus 5.0 and
OpenAI's GPT-5.6) were used subsequently to improve language and
readability, to draft and rewrite specific paragraphs under author
direction, and to assist in writing and debugging the project's Julia
code. All AI-generated output was reviewed and verified by the authors,
who take full responsibility for the final content.

\bibliographystyle{unsrtnat}
\bibliography{refs}

\end{document}